\documentclass[a4paper]{article}

\usepackage{RR}
\usepackage{amsfonts}
\usepackage{amssymb}
\usepackage{graphicx}
\usepackage{natbib}

\newtheorem{defn}{Definition}
\newtheorem{prop}{Proposition}
\newtheorem{rem}{Remark}
\newtheorem{hypo}{Hypotheses}
\newtheorem{thm}{Theorem}
\newtheorem{lem}{Lemma}

\def\qed{\relax\ifmmode\hskip2em \Box\else\unskip\nobreak\hskip1em $\Box$\fi}

\newcommand{\resp}{{\it resp. }}
\newcommand{\beq}{\begin{equation}}
\newcommand{\eeq}{\end{equation}}
\newcommand{\ie}{{\it i.e. }}

\newcommand{\eg}{{\it e.g. }}
\newcommand{\bi}{\begin{itemize}}
\newcommand{\ei}{\end{itemize}}
\newcommand{\dst}{\displaystyle}

\RRdate{Juillet 2006}
\RRauthor{
Ludovic Mailleret\thanks[inra]{URIH, INRA, 400 route des Chappes, BP\ 167, 06\ 903
   Sophia-Antipolis Cedex, France}
\and
Jean-Luc Gouz\'e\thanks[inria]{COMORE, INRIA, 2004 route des Lucioles, BP\ 93, 06\ 902,  Sophia-Antipolis Cedex, France}
\and
Olivier Bernard\thanksref{inria}}

\authorhead{Mailleret, Gouz\'e \& Bernard}

\RRtitle{Stabilisation globale d'une classe de syst\`emes positifs partiellement connus}
\RRetitle{Global Stabilization of a Class of Partially Known Positive Systems}

\RRresume{Nous consid\'erons ici un probl\`eme de stabilisation globale par retour de sortie d'une classe de syst\`emes non-lin\'eaires positifs de dimension \textit{n} poss\'edant une partie de forme analytique inconnue, mais cependant mesur\'ee. Nous proposons tout d'abord notre r\'esultat principal, un contr\^ole par retour de sortie, qui stabilise globallement le syst\`eme vers un \'equilibre r\'eglable appartenant \`a l'int\'erieur de l'orthant positif de l'espace d'\'etat. Bien qu'assez g\'en\'eral, ce r\'esultat est bas\'e sur des hypoth\`eses qui peuvent se r\'ev\'eler difficiles \`a v\'erifier sur un exemple concret. Dans un second temps, par l'interm\'ediaire d'un th\'eor\`eme sur une classe de syst\`emes positifs liant existence et stabilit\'e globale d'un \'equilibre fortement positif, nous proposons d'autres hypoth\`eses permettant d'appliquer notre r\'esultat. Ces nouvelles  hypoth\`eses, bien que plus restrictives, sont nettement plus faciles \`a v\'erifier. Le rapport s'ach\`eve par des exemples qui illustrent \`a la fois la potentielle complexit\'e des comportements en boucle ouverte qui peuvent \^etre produits par les syst\`emes consid\'er\'es (multi-stabilit\'e, cycle limite, chaos) et l'int\'er\^et de la proc\'edure de contr\^ole propos\'ee.}

\RRabstract{In this report we deal with the problem of global output feedback stabilization of a class of  $n$-dimensional nonlinear positive  systems possessing a one-dimensional unknown, though measured, part. We first propose our main result, an output feedback control procedure, taking advantage of measurements of the uncertain part, able to globally stabilize the system towards an adjustable equilibrium point in the interior of the positive orthant. Though quite general, this result is based on hypotheses that might be difficult to check in practice. Then in a second step, through a Theorem on a class of positive systems linking the existence of a strongly positive equillibrium to its global asymptotic stability, we propose other hypotheses for our main result to hold. These new hypotheses are more restrictive but much simpler to check. Some illustrative examples, highlighting both the potential complex open loop dynamics (multi-stability, limit cycle, chaos) of the considered systems and the interest of the control procedure, conclude this report.}

\RRmotcle{Syst\`emes positifs, syst\`emes partiellement connus, stabilisation globale.}

\RRkeyword{Positive systems, Partially know systems, Global stabilization.}

\RRprojet{Comore}

\RRtheme{\THBio}

\URSophia

\begin{document}
\makeRR

\section{Introduction}

Positive systems of ordinary differential equations are systems producing trajectories, that, if initiated in the positive orthant of $\mathbb{R}^n$, remain in this orthant for all positive time. Such systems are being payed much attention since they abound in many applied areas such as life sciences, social sciences, chemical sciences, telecommunications, traffic flows {\it etc}...  \citep{Luenberger1979,Farina2000,posta2003}.

In this paper we deal with the problem of global output feedback stabilization of a class of  $n$-dimensional nonlinear positive  systems possessing a 1-dimensional unknown, though measured, part. Our main result generalizes an earlier one obtained in lower dimensions and for simpler structures in the context of bioprocess control \citep{Mailleret2004a}.

This paper is organized as follows. We first introduce some  notations and definitions related to positivity for vectors and dynamical systems. We then introduce the class of uncertain positive systems that is concerned with our control problem.  We then propose an output feedback procedure, taking advantage of measurements of the uncertain part, able to globally stabilize the system towards an adjustable equilibrium point in the interior of the positive orthant. The hypotheses on which our main result is built might however be difficult to check in practice. In order to derive some hypotheses that are much simpler to check, we propose a theorem on a class of positive systems (the strongly positive concave cooperative systems) that links the existence of a positive equilibrium to its global asymptotic stability (GAS). Some examples found in the litterature illustrate the potential complex open loop dynamics that might be produced by the systems under study (multi-stability, cyclic behavior or even chaos). Simulations on periodic and chaotic behaviors stabilization illustrate the control efficiency.

\section{Preliminaries}

\subsection{Notations}

We consider the following autonomous nonlinear dynamical systems in $\mathbb{R}^n$: 
\begin{equation}\label{dyn_sys}
\dot x=f(x).
\end{equation}

For the sake of simplicity, we assume that the  functions  encountered throughout the paper are sufficiently smooth.

In the sequel, we will use the following notations:
\bi\item $x(t,x_0)$ denotes the forward orbit at time $t$ of system (\ref{dyn_sys}) initiated at $x_0$. 
\item $Df(x)$ denotes the Jacobian matrix of system (\ref{dyn_sys}) at state $x$.
\item $\mathcal{L}_{f}$ denotes the Lie derivative operator along the vector field defined by system (\ref{dyn_sys}).
\item $x_{|x_i=0}$ denotes a state vector whose $i-$th component equals zero.\ei

Moreover, a matrix is said to be {\it Hurwitz} if all of its eigenvalues have negative real parts.

\subsection{Positivity for vectors and dynamical systems} 

As we deal with positivity in $\mathbb{R}^n$ and some related concepts, we first clearly state what we refer to as positivity (\resp strong positivity) for vectors and dynamical systems. 

\begin{defn}[Positivity\label{pos_vec}]$~$\\
$x\in\mathbb{R}^n$ is positive (\resp strongly positive) and denoted $x\geq 0$ (resp.  $x\gg 0$) iff all its components are non-negative (\resp positive).\end{defn}

We denote $\mathbb{R}^n_+$ and $\textrm{int}(\mathbb{R}^n_+)$ the sets of positive vectors and strongly positive vectors belonging to $\mathbb{R}^n$, respectively. Same definitions of the relations $\geq$ and $\gg$ ({\it i.e.} component by component) will be used for matrices too. 

We now state the definition of positive systems of ordinary differential equations as well as the definition of a special class of positive systems, the strongly positive ones.

\begin{defn}[Positive Systems\label{pos_sys}\label{def_pos_sys}]$~$\\
System (\ref{dyn_sys}) is a positive system (\resp strongly positive) iff:
\begin{eqnarray*}
\forall i\in[1..n],~\forall x_{|x_i=0}\geq 0,~\dot{x_i}(x_{|x_i=0})=f_i(x_{|x_i=0})\geq
0,~
(resp.\ f_i(x_{|x_i=0})> 0).
\end{eqnarray*}\end{defn}

It is straightforward that Definition \ref{pos_sys} implies the following:

\begin{prop}\label{prop_pos_sys}$~$\\
Consider a positive system (\ref{dyn_sys}) (\resp strongly positive). Then:
\begin{eqnarray*}\forall (x_0,t)\in\mathbb{R}^n_+\times \mathbb{R}^+ ~~~(\resp~ \forall (x_0,t)\in  \mathbb{R}^n_+\times \mathbb{R}^+_*),~~~x(t,x_0)\geq 0 ~~(\resp \gg 0).
\end{eqnarray*}
\end{prop}

\begin{rem}
Though we have equivalence between Definition \ref{def_pos_sys} and Proposition \ref{prop_pos_sys} in the case of positive systems (this is the original definition of positive systems by \cite{Luenberger1979}), we only have an implication from Definition \ref{def_pos_sys} to Proposition \ref{prop_pos_sys} for strongly positive systems. Indeed systems with $f_i(x_{|x_i=0})\geq 0$ and a strongly connected Jacobian Matrix may produce strongly positive trajectories too \citep{Farina2000}. We however keep Definition  \ref{def_pos_sys} as the definition of strongly positive systems since this criterion is much easier to check on dynamical systems.
\end{rem}

\section{Main result\label{sec_main}}

\subsection{The considered class of  systems}

In the sequel we focus on the (global) asymptotic stabilization problem of  the following class of nonlinear positive uncertain systems in $\mathbb{R}^n$:
\begin{equation}\label{sys_gene}
\begin{array}{l}
\dot x = u f(x)+c \psi(x),\\
y=\psi(x).\end{array}
\end{equation}

\begin{itemize}
\item $u\geq 0$ being the scalar input of system (\ref{sys_gene}).
\item $f:~\mathbb{R}^n\rightarrow \mathbb{R}^n$
\item $c\in\mathbb{R}^n$
\item $\psi:~\mathbb{R}^n\rightarrow\mathbb{R}$ is an output of system (\ref{sys_gene})\end{itemize}
Function $\psi(.)$ features the {\it uncertain/unknown} part of the system. For control purposes, we suppose that, though analytically  unknown, $\psi(.)$ might be online measured. We moreover assume the following:
\begin{hypo}[H1]\label{hypo_sys_gene}\  
\begin{enumerate}
\item $\forall x\in\emph{\textrm{int}}(\mathbb{R}^n_+),\ \psi(x)> 0~$ and \ $\forall i ,~ c_i \psi(x_{|x_i=0})\geq 0$
\item $f(.)$ is such that the system $~\dot x=f(x)$ is positive
\item $\exists\beta_m\in \mathbb{R}^+$ such that $\forall
  \beta>\beta_m$, the system:\\[.15cm]
$\dot{x}=\beta f(x)+c.$\\[.15cm]
is a strongly positive system and possesses an equilibrium  $x^\star_\beta\gg 0$ which is GAS on $\mathbb{R}^n_+$.
\end{enumerate}\end{hypo}

These hypotheses will be commented after the statement of our main result. For now, just notice that Hypotheses H1-1 and H1-2 together with $\ u\geq 0\ $ ensure the positivity of the considered class of systems  (\ref{sys_gene}).

\subsection{Main result}

\begin{thm}\label{theo_main_result}$~$\\ Consider a system (\ref{sys_gene}) and let Hypotheses H1 hold. Then for all  $\gamma>\beta_m$, the nonlinear control law:
\beq\label{control_law} u=\gamma y=\gamma \psi(x).
\eeq
globally stabilizes system (\ref{sys_gene}) on $\emph{\textrm{int}}(\mathbb{R}^n_+)$ towards the strongly positive equilibrium $x^\star_\gamma$ such that:
\[f(x^\star_\gamma)=\dst\frac{-1}{\gamma}\ c.
\]
\end{thm}

To prove Theorem \ref{theo_main_result}, we need  the following Lemma on Lyapunov functions for  positive systems possessing a  strongly positive GAS equilibrium.

\begin{lem}\label{lyapu}$~$\\
A positive system (\ref{dyn_sys}) possesses a strongly positive equilibrium point $x^\star$ GAS on $\emph{\textrm{int}}(\mathbb{R}^n_+)$ iff there exists a smooth real valued function $V(x)$  such that:
\[\begin{array}{l}
\forall x\gg 0, x\neq x^\star, V(x)>0,\\
V(x^\star)=0,\\
\forall x\gg 0, x\neq x^\star, \dot V(x)<0,\\
\dst\bigcup_{\alpha>0}\{x,~V(x)\leq\alpha\}=\emph{\textrm{int}}(\mathbb{R}^n_+).
  \end{array}
\]\end{lem}

{\bf Proof}{\bf $~$(Lemma \ref{lyapu})}\\
Consider a positive system (\ref{dyn_sys}) and suppose it possesses a strongly positive equilibrium $x^\star$ GAS on $\textrm{int}(\mathbb{R}^n_+)$.\\
Now consider the following increasing (in the $\gg$ sense) change of variables:
\[z=\ln\left(\dst\frac{x}{x^\star}\right)=
\left(\begin{array}{c}
\ln\left(\dst\frac{x_1}{x_1^\star}\right)\\
\vdots\\
\ln\left(\dst\frac{x_n}{x_n^\star}\right)
      \end{array}\right).
\]
that maps the strongly positive orthant on $\mathbb{R}^n$ and $x^\star$ to 0. Then we have:
\begin{equation}\label{chg_var_z}
\dot z=\dst\frac{f(x^\star e^z)}{x^\star}e^{-z}.
\end{equation}
Since $x^\star$ is a GAS equilibrium on $\textrm{int}(\mathbb{R}^n_+)$ for (\ref{dyn_sys}), so is 0 on $\mathbb{R}^n$ for system (\ref{chg_var_z}). Then, from Kurzweil's converse Lyapunov theorem on GAStability (see Theorem 7 in  \citep{Kurzweil1956}), there exists a smooth radially unbounded Lyapunov function $W(z)$ for system (\ref{chg_var_z}). Coming back to the original variables, consider the function:
\[V(x)=W(\ln(x)-\ln(x^\star)).\]
Since $W(z)$ is a Lyapunov function for system (\ref{chg_var_z}), we have:
\[\begin{array}{l}
\forall x\in\textrm{int}(\mathbb{R}^n_+),~x\neq x^\star,~V(x)=W(z)>0,\\
V(x^\star)=W(0)=0,\\
\forall x\in\textrm{int}(\mathbb{R}^n_+),~x\neq x^\star,~\dot V(x)=\dot W(z)<0.\end{array}\]
Moreover, since $W(z)$ is smooth so does $V(x)$ and since $W(z)$ is  radially unbounded on $\mathbb{R}^n$, we have:
\[\bigcup_{\beta>0}\{z,\ W(z)\leq \beta\}=\mathbb{R}^n.\]
Thus:
\[\bigcup_{\alpha>0}\{x,\ V(x)\leq \alpha\}=\textrm{int}(\mathbb{R}^n_+).\]
which concludes the first part of the proof.

The reverse implication is easily obtained through classical Lyapunov theory. \qed

{\bf Proof}{\bf $~$(Theorem \ref{theo_main_result})}\\ 
Consider the positive system (\ref{sys_gene}) under control law (\ref{control_law}). We get:
\begin{equation}\label{closed_loop_sys}
\dot x =\psi(x)(\gamma f(x)+c) \triangleq g(x).
\end{equation}
Since $\gamma>\beta_m$ we know from H1-3 that the system:
\begin{equation}\label{time_chg_closed_loop}\dot x=\gamma f(x)+c \triangleq h(x).
\end{equation}
is strongly positive and possesses a (strongly positive) GAS equilibrium $x^\star_\gamma$. Then, from Lemma \ref{lyapu}, for all $\gamma>\beta_m$, system (\ref{time_chg_closed_loop}) admits a real valued function $V_\gamma(x)$ verifying Lemma \ref{lyapu}  properties.

Now we consider system (\ref{closed_loop_sys}). From H1-1, it is clear that we have:
\[\begin{array}{l}
\forall x\gg 0, x\neq x^\star_\gamma,\ V_\gamma(x)>0,\\
V_\gamma(x^\star_\gamma)=0.\end{array}\]
Moreover:
\[\mathcal{L}_g V_\gamma(x)=\psi(x)\mathcal{L}_h V_\gamma(x).
\]
so that since H1-1 holds: 
\[\forall x\gg 0, x\neq x^\star_\gamma,\ \mathcal{L}_{g} V_\gamma(x)<0.
\]
and the following still holds:
\[\dst\bigcup_{\alpha>0}\{x,\ V_\gamma(x)\leq\alpha\}=\textrm{int}(\mathbb{R}^n_+).
\]
So that we can apply Lemma \ref{lyapu} to the controlled system (\ref{closed_loop_sys}). Then the equilibrium $x^\star_\gamma$ is GAS for the controlled system (\ref{closed_loop_sys}) on  $\textrm{int}(\mathbb{R}^n_+)$. \qed

\begin{rem}
It is important to notice that the proposed control law does not, to be applied, require any {\it a priori} analytical or quantitative knowledge of the function $\psi(x)$. The only two requirements are (i) its positivity (a qualitative  property), and (ii) the possibility to measure this quantity on-line. This is particularly important for instance in biological process control, when some parts of the model are only qualitatively known (usually in terms of sign) while some other parts are precisely (analytically) known. The special case of bioprocesses presented by \cite{Mailleret2004a} does illustrate this point and shows the interest of such control procedures for real life bioprocesses management. We come back to this point later in the examples section. \end{rem}

\subsection{Comments on Hypotheses H1}

As we have previously noted, Hypotheses H1-1 and H1-2 guarantees the positivity of the systems (\ref{sys_gene}).

First part of Hypothesis H1-1 defines a qualitative property of the unknown part of the system, which, though one-dimensional, possibly acts on all state variables dynamics via the vector $c$. Notice that this is a really loose hypothesis since the only required property is $\psi(x)$ positivity as $x$ is strongly positive. No specific analytical form is assumed on $\psi(.)$. The second part of H1-1 is required to ensure system (\ref{sys_gene})'s positivity as the input $u$ equals 0. H1-2 is required for the system's positivity as the input $u$ is positive. 

Hypothesis H1-3 is a stronger one. It reads that if the uncertain part $\psi(.)$ were positive and constant, then there must exist an input infimum $\beta_m$ above which system (\ref{sys_gene}) would be strongly positive and would possess a (strongly positive) GAS equilibrium. This property might recall ``minimum phase" 
conditions, frequenlty encountered when dealing with the stabilization of nonlinear and/or uncertain systems. Minimum phase conditions  ensures the GAStability of an equilibrium as the system's outputs are constant. Minimum phase based control results usually stabilises the output(s) in order to conclude on the whole system, using results on cascaded nonlinear systems' stability. However uncertainties for these systems are in general not on the outputs' dynamics, but rather in the intrinsic ``stable" part. A ``minimum phase" approach can thus not be used in our case since the output $\psi(x)$ and its dynamics is, to a large extent, unknown.

From the analytical expression of the system that is to be stabilized, both Hypotheses H1-1 and H1-2 are easy to check. This is however not the case for H1-3. The existence of a unique and  strongly positive equilibrium of:
\[\dot x = \beta f(x)+c\]
might indeed be easily proven, but the demonstration of its global stability is in general a much harder task. Nevertheless, for some special classes of positive systems, one can link the existence of a strongly positive equilibrium to its global asymptotic stability.  This is for instance the case for positive linear systems, see \eg  \citep{Luenberger1979,Farina2000}. It is clear that such a property is of great interest: it allows to prove in one step both the existence of a strongly positive equilibrium and its GAStability \ie Hypothesis H1-3. Some works can be found in the litterature regarding this problem for nonlinear systems: see \citep{Silva-Navarro1997,De2001a}. These results can indeed be used to check H1-3. It is to be noticed that these works are all related to cooperative systems (see Definition \ref{cosy} in the sequel). 

In the next section we also base our approach on cooperative systems. We propose a class of positive systems for which the existence of a strongly positive equilibrium implies its GAStability. We then derive a set of ``easy to check" hypotheses  that are sufficient for our control procedure (\ref{control_law}) to be applied.

\section{Sufficient conditions to verify H1}

\subsection{Useful known results}

We now recall the Metzler matrices definition introduced by \cite{Metzler1945} in mathematical economy.

\begin{defn}[Metzler Matrix]$~$\\
A matrix is Metzler iff all its off-diagonal elements are non-negative.
\end {defn}

We now focus on two results on Metzler matrices that are consequences of the original theorem of Perron-Frobenius on positive matrices \citep{Frobenius1908}. The first comes from \cite{Luenberger1979}: 
\begin{thm}\label{metz_luen}$~$\\
Consider a Metzler matrix $A$ and a vector $b\in\emph{\textrm{int}}(\mathbb{R}^n_+)$; Then, $A$ is Hurwitz  iff: 
\[\exists x\geq 0,~Ax+b=0.\]
\end{thm}

The second can be found in \citep{Smith1986b}:
\begin{thm}\label{metz_smith}$~$\\
Consider a Metzler matrix $A$. Then: 
\begin{itemize}
\item The dominant eigenvalue of $A$ (\ie of largest real part) is real; its associated eigenvector is positive.
\item consider a matrix $B\geq A$. Then the real part of the dominant eigenvalue of $B$ is greater or equal to the dominant eigenvalue of $A$.  
\end{itemize}\end{thm}

Let us recall the definition of the ``cooperative systems", introduced by \cite{Hirsch1982} (see also \citep{Smith1995}) and that will be  payed attention in the sequel. 
\begin{defn}[Cooperative Systems]\label{cosy}$~$\\ System (\ref{dyn_sys}) is a cooperative system iff its Jacobian $Df(x)$ is Metzler for all $x\in\mathbb{R}^n$.
\end{defn}

One of the most interesting result on cooperative systems is to be found in \citep{Smith1995} and reads: 
\begin{thm}\label{order}$~$\\
Consider a cooperative system (\ref{dyn_sys}) and two initial conditions $x_0$ and $y_0$ in $\mathbb{R}^n$ such that $x_0\geq y_0$ (\resp $\gg$). Then:
\[\forall t\geq 0,~x(t,x_0)\geq y(t,y_0)~~(\resp \gg).\]
\end{thm}

\subsection{Preliminary result}

\begin{prop}[Positive Cooperative Systems] \label{sys_pos_coop}$~$\\ 
Consider a cooperative system (\ref{dyn_sys}). Then, the following are equivalent:
\bi\item system (\ref{dyn_sys}) is positive  (\resp strongly positive)
\item $f(0)\geq 0$ (resp. $f(0)\gg 0$)\ei
\end{prop}

{\bf Proof}\\ Suppose $f(0)\geq 0$ (resp. $f(0)\gg 0$). Let us denote the $i$-th row of $Df(x)$ by $Df_i(x)$. Then we have for all $x_{|x_i=0}\geq 0$: 
\begin{eqnarray}\label{fond_calc}
f_i(x_{|x_i=0})=
f_i(0)+\left[\ \int_0^1Df_i(sx_{|x_i=0})ds\right].\ x_{|x_i=0}.\end{eqnarray}
Since $Df(x)$ is Metzler for all $x$, the only possible negative term in the scalar product in (\ref{fond_calc}) is zero: indeed, $Df_{i,i}$ (the sole possible negative term of $Df_i$) is multiplied by $x_i=0$. Then: 
\[\forall x_{|x_i=0}\geq 0,~f_i(x_{|x_i=0})\geq f_i(0)\geq 0~~ (resp. > 0).\]
Since this holds for all $i\in[1..n]$, we conclude, using Proposition \ref{def_pos_sys}, that the system is positive ({\it resp.} strongly positive).

The reverse implication is obvious.\qed%\end{pf}

\begin{rem} Notice that the class of positive cooperative systems is a generalization of the class of linear positive systems  \citep{Luenberger1979,Farina2000}.

\end{rem}

We now state a result on global asymptotic stability of an equilibrium point for a special class of strongly positive cooperative systems.
\begin{thm}\label{theo_sys_coop1}$~$\\
Consider a strongly positive concave cooperative (SPCC) system   (\ref{dyn_sys}), \ie that verifies the following condition:
\begin{equation}\label{concavity}
\forall (x,y)\in \mathbb{R}^n_+\times\mathbb{R}^n_+,~x \leq y \Rightarrow Df(x)\geq Df(y).\end{equation}
Then, if system (\ref{dyn_sys}) has a positive equilibrium $x^\star$, it is single, strongly positive and GAS on $\mathbb{R}^n_+$.
\end{thm}

{\bf Proof:}\\ We first show the exponential stability of a positive equilibrium $x^\star$ (that is necessarily strongly positive) of a SPCC system (\ref{dyn_sys}). 

We obviously have:
\beq\label{zoubida1} f(x^\star)=0=f(0)+\left[\ \int_0^1 Df(sx^\star)ds\right]x^\star.\eeq 
Since $Df(x)$ is Metzler for all positive $x$, so does the  bracketed matrix (denoted $F(x^\star)$ in the sequel) in equation (\ref{zoubida1}). Since the considered system is a strongly positive cooperative system, we have from Proposition \ref{sys_pos_coop}: $f(0)\gg 0$. Applying Theorem \ref{metz_luen}, we conclude that matrix $F(x^\star)$ is Hurwitz. On the other hand,  equation (\ref{concavity}) implies:
\[\forall s\in [0,~1],~Df(sx^\star)\geq Df(x^\star).\]
which yields:
\[F(x^\star)=\left[\ \int_0^1 Df(sx^\star)ds\right]\geq Df(x^\star).\]
Since $F(x^\star)$ is a Hurwitz Metzler matrix, we show using Theorem \ref{metz_smith} that its dominant eigenvalue is negative real. From the above inequation and Theorem \ref{metz_smith}, we conclude that $Df(x^\star)$'s dominant eigenvalue is negative real too, proving that $x^\star$ is exponentially stable.

We now show the unicity of a strongly positive equilibrium $x^\star$. 

Suppose that the considered system possesses two strongly positive
equilibria $x_1^\star$ and $x_2^\star$, such that $x_1^\star\geq x_2^\star$. Then: 
\beq\label{2eqlies}
\begin{array}{ll}
0&=f(x_1^\star)-f(x_2^\star),\\[.1cm]
&=
\displaystyle{\left[\ \int_0^1Df(sx_1^\star+(1-s)x_2^\star)ds\right](x_1^\star-x_2^\star)}.\end{array}\eeq
From equation (\ref{concavity}) and since $x_1^\star\geq x_2^\star$, we use  similar arguments than in the previous part of the proof to show that:  
\[\left[\ \int_0^1Df(sx_1^\star+(1-s)x_2^\star)ds\right]\leq Df(x_2^\star). \]
Since $x_2^\star$ is a strongly positive equilibrium, $Df(x_2^\star)$ is Hurwitz. Then, using Theorem \ref{metz_smith}, the bracketed matrix is Hurwitz too. Then it possesses an inverse in $\mathcal{M}_n(\mathbb{C})$ which, together with equation (\ref{2eqlies}), implies that $x_1^\star=x_2^\star$. 

Suppose now that the considered system possesses two strongly positive equilibria $x_1^\star$ and $x_2^\star$, that are not linked by the relation ``$\geq$'', \ie ($x_1^\star-x_2^\star)$ is neither positive, nor negative. Let us define the parallelotopes $\mathcal{B}_z\subset\mathbb{R}^n$, such that:
\[\forall z \in \mathbb{R}^n_+,~\mathcal{B}_z=\left\{x\in\mathbb{R}^n_+,~ x\leq z\right\}\] 
Using Theorem \ref{order} together with the considered system's  cooperativity and strong positivity, it is easy to show that $\mathcal{B}_{x_1^\star}$ and $\mathcal{B}_{x_2^\star}$ are positively invariant sets. We now define the state $x_3$, such that:
\[x_{3,i}=\min(x_{1,i}^\star~,~x_{2,i}^\star)\]
Since 
$\mathcal{B}_{x_3}=\mathcal{B}_{x_1^\star}\cap\mathcal{B}_{x_2^\star}$, it is positively invariant by system (\ref{dyn_sys}) too. Then, using Brouwer fixed point Theorem (see \eg \cite{Zeidler1985}), there must exist a third equilibrium $x_3^\star$ belonging to $\mathcal{B}_{x_3}$. Note that, from $x_3$ definition, we necessarily have:
\[x_3^\star\leq x_1^\star~\textrm{ and }~x_3^\star\leq x_2^\star\]
This case has been treated previously and we conclude that $x_1^\star=x_3^\star=x_2^\star$, what yields a contradiction.

We conclude that if a SPCC system (\ref{dyn_sys}) possesses a positive equilibrium, it is the sole positive equilibrium.

We achieve the proof of Theorem \ref{theo_sys_coop1} showing the global attractivity  (on $\mathbb{R}^n_+$) of the positive equilibrum $x^\star$ of a SPCC system (\ref{dyn_sys}).

Consider the following real valued function: 
\[\begin{array}{llll}
\forall i\in[1..n],~g_i: &\mathbb{R}^+&\rightarrow &\mathbb{R}\\
& k&\mapsto &f_i(kx^\star)\end{array}\]
Since $f(0)\gg0$ and $f(x^\star)=0$, it is clear that $g_i(0)>0$ and $g_i(1)=0$. Then, there exists a $k_0\in(0,~1)$ such that: \[\frac{dg_i}{dk}(k_0)=g_i(1)-g_i(0)<0.\]
From equation (\ref{concavity}), we get:
\[0\leq k_0\leq k\Rightarrow 
0>\frac{dg_i}{dk}(k_0)\geq\frac{dg_i}{dk}(k).\]
showing that function $g_i$ strictly decreases from $k_0$ to $+\infty$. Since $g_i(1)=0$, $g_i(k)<0$ for all $k>1$. This holds for all $i\in[1..n]$, then: 
\[\forall k>1,~f(kx^\star)\ll 0\]
Untill the end of the proof, we assume $k>1$. Notice that, from the cooperativity and strong positivity properties, 
$\mathcal{B}_{kx^\star}$ is a positively invariant set. Now consider $\ddot x$, the time derivative of $\dot{x}$. We have: 
\[\ddot x=Df(x)\dot x\]
Since $Df(x)$ is a Metzler matrix for all $x$, it is straightforward that if the initial state velocity $\dot x(t=0)$ is positive (\resp negative) then it will remain positive (\resp negative) for all positive time.

Remind that $\dot x(x=0)\gg 0$ and $\dot x(x=kx^\star)\ll 0$. 
Then, the forward trajectory of system (\ref{dyn_sys}) initiated at $x_0=0$ (\resp $x_0=kx^\star$) is, in the ``$\geq$'' sense, an increasing (\resp decreasing) function of the time. Moreover, from $\mathcal{B}(kx^\star)$ positive invariance, it is lower (\resp upper) bounded by $0$ (\resp $kx^\star$), thus it converges, necessarily toward an equilibrium. 

$x^\star$ is the sole equilibrium belonging to $\mathbb{R}^n_+$, the two considered trajectories converge thus toward $x^\star$. From system (\ref{dyn_sys}) cooperativity together with Theorem \ref{order}, we conclude that each trajectory initiated in $\mathcal{B}_{kx^\star}$ converges toward $x^\star$. Since this holds for any $\mathcal{B}_{kx^\star}$ ($k>1$), we have shown that $x^\star$ is globally attractive on $\mathbb{R}^n_+$, which concludes the proof.\qed

\begin{rem}
Proof of theorem \ref{theo_sys_coop1} uses to a large extent the ideas proposed by \cite{Smith1986b}. However, his proof was dedicated to another class of systems:  Kolmogoroff-type population dynamics models of cooperating species, which are based on different hypotheses. One can also have a look at \cite{Smith1986a} for  some results related to Theorem \ref{theo_sys_coop1}.\end{rem} 

\begin{rem} It is clear that Theorem \ref{theo_sys_coop1} is quite similar to Hirsch's results \citep{Hirsch1985} on cooperative systems, rephrased by \cite{Piccardi2002}, and stating that for a bounded strongly monotone systems(\eg cooperative systems verifying $x_0\geq y_0 \Rightarrow x(t,x_0)\gg y(t,y_0)$ , the trajectories, apart from a set of zero measure (thus at most $n-1$ dimensional) of initial conditions, converge to the set of equilibrium points. Thus if such a system were bounded and the equilibrium unique then almost all trajectories converge to it. Notice also that Smith  provides a related Theorem (see \citep{Smith1995}, theorem 3.1). 
Here for SPCC systems, both uniqueness of the equilibrium (if it exists) and (real) global convergence are guaranteed, which is of particular interest for the application of our control procedure (\ref{control_law}).
\end{rem}

\subsection{Sufficient conditions to verify H1}

On the basis of SPCC systems and Theorem \ref{theo_sys_coop1} we provide sufficient conditions that are much simpler to check and that guarantees that Hypotheses H1 are satisfied.

Consider a system (\ref{sys_gene}) under the following hypotheses:
\begin{hypo}[H2]\label{hypo_sys_gene2}\  
\begin{enumerate}
\item $\forall x\in\emph{\textrm{int}}(\mathbb{R}^n_+),\ \psi(x)> 0~$ and  $\forall i ,~ c_i \psi(x_{|x_i=0})\geq 0$ 
\item $f(0)\geq 0$ 
\item $f(.)$ is such that system $~\dot x=f(x)$ is cooperative
\item $\forall x_1, x_2 \in \mathbb{R}^n_+, ~x_1\leq x_2 \Rightarrow Df(x_1)\geq Df(x_2)$
\item $\exists\beta_m\in \mathbb{R}^+$ such that $\forall
  \beta>\beta_m,~\beta f(0)+c\gg 0$
\item $\forall \beta>\beta_m,~\exists x^\star_\beta \in\emph{\textrm{int}}(\mathbb{R}^n_+)$ such that $\beta f(x^\star_\beta)+c= 0$
\end{enumerate}\end{hypo}

One can show that H2 implies H1, thus that a system (\ref{sys_gene}) under H2 is a candidate for Theorem \ref{theo_main_result}.

Indeed, H2-1 and H1-1 are the same. H2-2 together with H2-3 implies through Proposition \ref{sys_pos_coop} that $\dot x=f(x)$ is a positive system \ie H1-2 holds. H2-3 together with H2-5 and Proposition \ref{sys_pos_coop} implies that for all $\beta>\beta_m$, the system:
\begin{equation}\label{remark2}
\dot x=\beta f(x)+c.\end{equation}
is a strongly positive cooperative system. This fact together with H2-4 implies that (\ref{remark2}) is a SPCC system. H2-6 ensures system \ref{remark2} has, for all $\beta>\beta_m$, an equilibrium $x^\star_\beta$ that is, through Theorem \ref{theo_sys_coop1}, GAS on $\mathbb{R}^n_+$ so that H1-3 holds.\qed

In the next section we provide numerical simulations, based on examples fulfilling hypotheses H2, that illustrates Theorem \ref{theo_main_result} interest and efficiency.

\section{Illustrative Examples\label{sec_example}}

We show on three examples that the possible open loop dynamics of  systems (\ref{sys_gene}) under H2 might be very complex. Proofs of the complex dynamics are not detailed as they do not fall beyond the scope of this paper, but they can easily be performed with the help of the proposed references.

{\bf Example 1} (Bi-stability)\label{example_1} This system models a 
simple bioreaction occuring in a continuous stirred tank reactor with substrate inhibitory effects (see {\it e.g.} \citep{Smith1995a}). This is a simple model of an anaerobic digester, an apparatus used for waste water treatment \citep{Mailleret2004a}. Both variables represent concentrations inside the reactor: $x_1$ denotes the substrate (pollutant) consumed by the biomass (anaerobic microorganisms) $x_2$ to grow at a per capita rate $\mu(x_2)$. $u$ denotes the ``dilution rate'' (\ie passing flow per volume unit) feeding the reactor with substrate at a concentration $x_{1,in}$ and withdrawing a blend of $x_1$ and $x_2$ from it. We get the following system that will be refered to as $(\mathcal{S}_1)$ in the sequel.
\begin{eqnarray*}\begin{array}{rl}
\dot x & =u\left[\left(\begin{array}{cc}
-1 & 0 \\ 0 & -1
\end{array}\right) x+\left(\begin{array}{c}
x_{1,in}\\ 0\end{array}\right)\right]+
\left(\begin{array}{c}
-k\\ 1\end{array}\right)\mu(x_1)x_2,\\[.5cm]
&\triangleq uf_1(x)+c_1\psi_1(x).\end{array}\end{eqnarray*}\vspace*{-.25cm}
with: $\mu(x_1)=\displaystyle\frac{\mu_m x_1}{K_m+x_1+x_1^2/K_i}$

$x_{1,in},k,\mu_m, K_m$ and $K_i$ are positive constants. 

It is clear that Hypotheses H2-1, H2-2, H2-3 and H2-4 hold. Moreoever H2-5 holds with $\beta_m=k/x_{1,in}$. Now pick a $\beta>\beta_m$, then $\beta f_1(x)+c_1\psi_1(x)=0$ 
possesses a solution $x^\star_\beta$ such that: 
\[\forall \beta>\beta_m,~x^\star_\beta=\left(\begin{array}{c}
x_{1,in}-\dst\frac{k}{\beta}\\
\dst\frac{1}{\beta}
  \end{array}\right)\gg 0.\]
so that H2-6 holds. Then we have through Theorem \ref{theo_main_result}:
\begin{prop} For all $\gamma>k/x_{1,in}$ the control law 
\beq\label{con_dig_ana} u(.)=\gamma \psi_1(x)\eeq
globally stabilizes, on \emph{int}$(\mathbb{R}^n_+)$, system $(\mathcal{S}_1)$ towards $x^\star_\gamma$, the (sole) solution of $\gamma f_1(x)+c_1=0$.
\end{prop}

\begin{figure}[!h]\begin{center}
\includegraphics[width=\columnwidth]{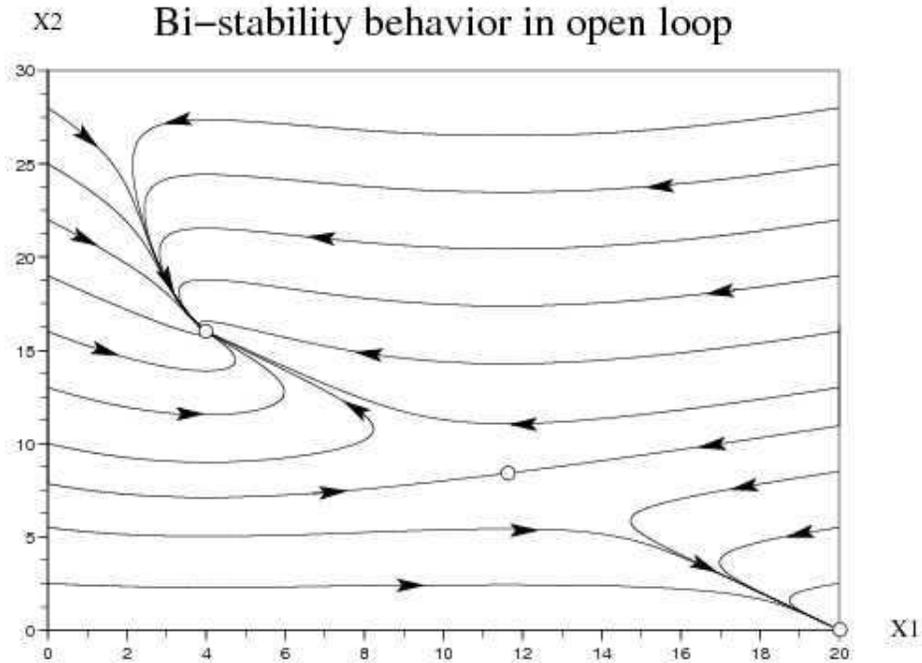}\end{center}  
\caption{State space representation of various open loop forward trajectories of $(\mathcal{S}_1)$ with \texttt{o} representing equilibria and  $u$ belonging to $(\mu(x_{1,in}), \max(\mu))$: existence of two locally stable equilibria. \label{fig_example_1}}  
\end{figure}

Now let us consider the open loop behavior that might be produced by $\mathcal{S}_1$. Suppose that $x_{1,in}>\textrm{argmax}\ \mu$ and $u\in (\mu(x_{1,in}),\max(\mu))$ is constant. Then the considered system  possesses three equilibria, two of which being stable on their respective basins of attraction
separated by the stable manifold of the third equilibrium that is a saddle point. This behavior is depicted on figure \ref{fig_example_1}; see \eg \citep{Mailleret2004} for a rigorous analysis.

The two stable equilibria correspond for the strongly positive one to the ``operating point" (the biomass survives inside the reactor and the pollutant concentration is reduced compared to its in the inflow) while the other do not (the biomass disappears and the waste water is no more treated). It is then necessary to ensure that no trajectory may be driven to this latter equilibrium. In this real life example, the quantity $\psi_1(x)$ can be easily measured since it corresponds to the outflow of biogaz that is produced by the anaerobic digester. Taking adavantage on this measurement, control procedure (\ref{control_law}) has been applied to a real life pilot scale anaerobic digester. It has shown its interest and efficiency for the global stabilization of such multistable devices. It thus avoids the crash of the system, see \citep{Mailleret2004a} for more details.

{\bf Example 2} (Attractive limit cycle\label{example_2})\\ The following  system is analogous to Goodwin metabolic models \cite{Goodwin1965}. It will be refered to as system  $(\mathcal{S}_2)$ in the sequel.
\begin{eqnarray*}\begin{array}{rl}
\dot x &=u\left(\begin{array}{c}
-lx_1\\
\frac{\mu_1 x_1}{k_1+x_1}-x_2+\alpha_1\\
\frac{\mu_2 x_2}{k_2+x_2}-x_3+\alpha_2
\end{array}
\right)
+\left(\begin{array}{c}
1\\0\\0
\end{array}\right)\dst\frac{1}{1+x_3^{n}},\\[.8cm]
& \triangleq uf_2(x)+c_2\psi_2(x).
\end{array}
\end{eqnarray*}
with $l=2.1$, $\mu_1=2/2.1$, $\mu_2=4*(0.01+1/2.1)$, 
$k_1=1/4.2$, $k_2=0.01+1/2.1$, $\alpha_1=0.01$, 
$\alpha_2=1-2(0.01+1/2.1)$  and $n=80$. It is easily proven that H2-1, H2-2, H2-3 and H2-4 hold true. Moreover 
H2-5 holds with $\beta_m=0$. Now pick a $\beta>\beta_m$, then consider the equation $\beta f_2(x)+c_2\psi_2(x)=0$. We get that $x_3$ must be a solution of:
\[x_3=\alpha_2+\dst\frac{\mu_2(\mu_1+\alpha_1(k_1\beta l(1+x_3^n)+1))}{\mu_1+(\alpha_1+k_2)(k_1\beta l(1+x_3^n)+1)}\]
evaluating both sides of this equation at $x_3=0$ and $+\infty$, we show that there exist a solution $x_{\beta,3}^*$ which is positive for all $\beta>\beta_m$. $x_{\beta, 2}^*$ and $x_{\beta,1}^*$ are then computed and shown to be positive too, so that H2-6 is  verified. 

\begin{figure}[!h]\begin{center}
\includegraphics[width=\columnwidth]{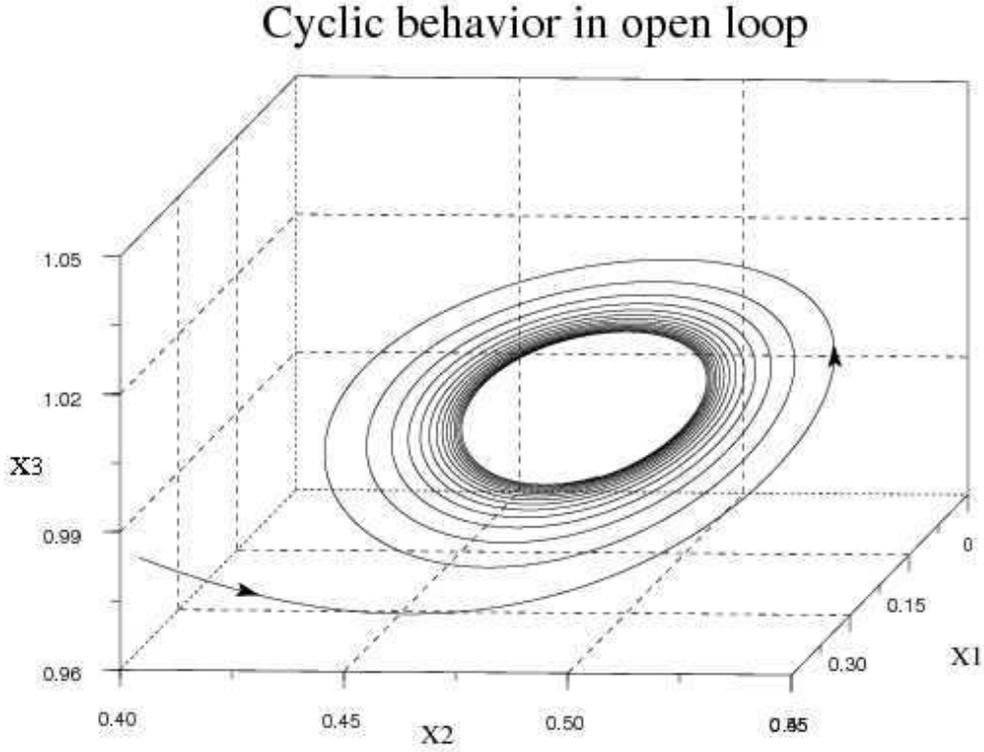}\end{center}
\caption{State space representation of a forward trajectory of $(\mathcal{S}_2)$ with $u=1$: existence of an attractive limit cycle.\label{fig_example_2}} 
\end{figure}

Suppose that the considered system operates in open loop with the non-negative input $u=1$. Then, as depicted on figure \ref{fig_example_2}, $(\mathcal{S}_2)$ possesses an attractive limit cycle around an unstable equilibrium point. The existence of a nontrivial periodic orbit can be proven with Theorem  1 from \citep{Hastings1977}.

\begin{figure}[!h]\begin{center}
\includegraphics[width=\columnwidth]{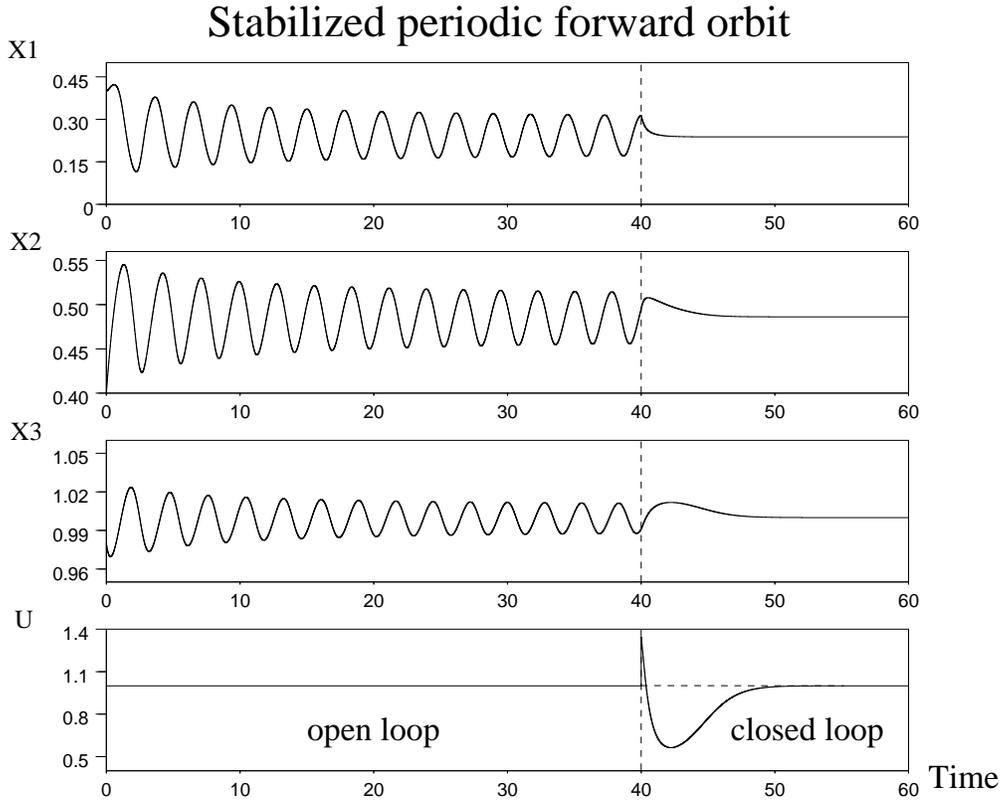}\end{center}
\caption{Time varying representation of a forward trajectory of $(\mathcal{S}_2)$ with $u=1$ before $t=40$ and under  control law (\ref{con_cyc_lic}) from $t=40$ to the end, with $\gamma=2>\beta_m(=0)$ computed s.t. asymptotically $u(.)=\gamma\psi_2(x)$ equals its open loop value.\label{fig_example_2_con}}
\end{figure}

As illustrated on Figure \ref{fig_example_2_con}, we get from {Theorem\ \ref{theo_main_result}:}
 \begin{prop} For all $\gamma>0$ the control law 
\beq\label{con_cyc_lic} u(.)=\gamma \psi_2(x)\eeq
globally stabilizes, on \emph{int}$(\mathbb{R}^n_+)$, system $(\mathcal{S}_2)$ towards $x^\star_\gamma$, the (sole) solution of $\gamma f_2(x)+c_2=0$.
\end{prop}

\vspace*{-.1cm}
{\bf Example 3} (Chaos\label{example_3}) The following system is a state space translation of the model of a chemical system so called ``auto-catalator'', see \citep{Petrov1992}. It reads:
\begin{eqnarray*}\begin{array}{rl}
\dot x\!&=\!u\!\left[\!\left(\begin{array}{ccc}
-1&0&k_2\\
\dst\frac{1}{k_1}&\dst\frac{-1}{k_1}&0\\
0&1&-1
\end{array}\!\right)\!x\!+\!
\left(\!\begin{array}{c}
k_2(k_3\!-\!k_4)\\
0\\
k_4
\end{array}\!\right)\!\right]%~~~\\
\!+\!\left(\!\begin{array}{c}
-1\\
\dst\frac{1}{k_1}\\
0
\end{array}\!\right)\!x_1x_2^2,\\[.8cm]
&\triangleq u f_3(x)+c_3\psi_3(x).\end{array}
\end{eqnarray*}
$k_1$ is positive, $k_2\in(0,\ 1)$ and $k_3>k_4>0$ . 

It is clear that H2-1, H2-2, H2-3 and H2-4 hold. Moreover 
H2-5 holds with $\beta_m=1/(k_2(k_3-k_4))$. Now pick a $\beta>\beta_m$, then $\beta f_3(x)+c_3\psi_3(x)=0$ 
possesses a solution $x^\star_\beta$ such that:
\[\forall \beta>\beta_m,~x^\star_\beta=\left(
\begin{array}{c}
\dst\frac{\beta k_2 k_3+k_2-1}{\beta(1-k_2)}\\[.25cm]
\dst\frac{k_2 k_3}{1-k_2}\\[.25cm]
\dst\frac{k_2(k_3-k_4)+k_4}{1-k_2}\end{array}
\right)\gg 0.\]
since $k_2\in(0,\ 1)$ and $k_3>k_4>0$. 
Then H2-6 holds true.

\begin{figure}[!h]\begin{center}
\includegraphics[width=\columnwidth]{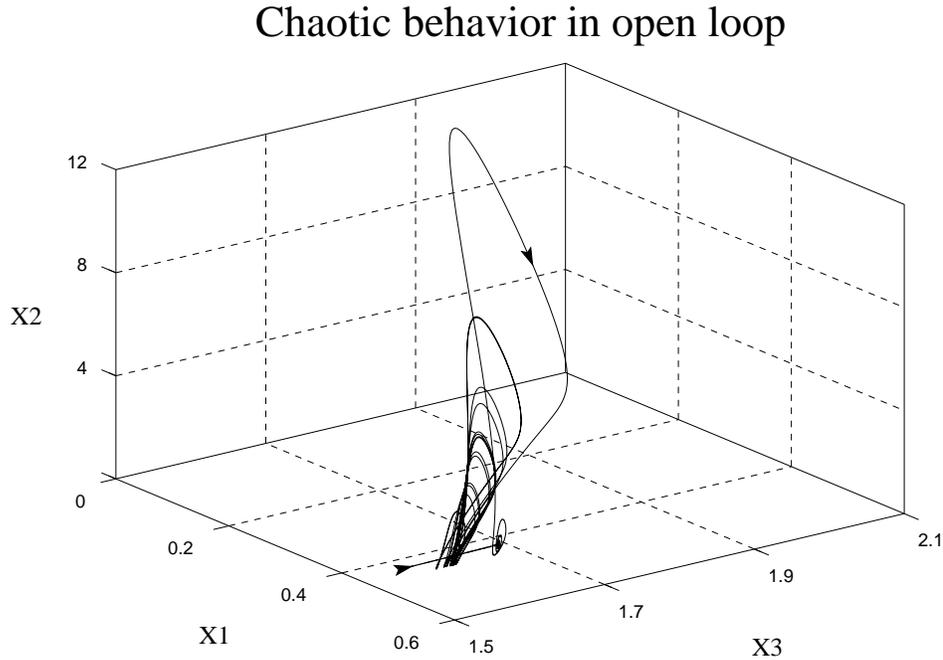}\end{center}
\caption{State space representation of a forward chaotic trajectory of the system $(\mathcal{S}_3)$ with $u=1$. \label{fig_example_3}}
\end{figure}

\cite{Petrov1992} prove that Example \ref{example_3} exhibits a chaotic behavior for the parameter values $k_1=0.015,~k_2=0.301,~k_3=2.5$ and $u=1$ for all $k_4$. As a special case, we choose $k_4=0.56$. An example of a forward chaotic trajectory is shown on figure \ref{fig_example_3}.

\begin{figure}[!h]\begin{center}
\includegraphics[width=\columnwidth]{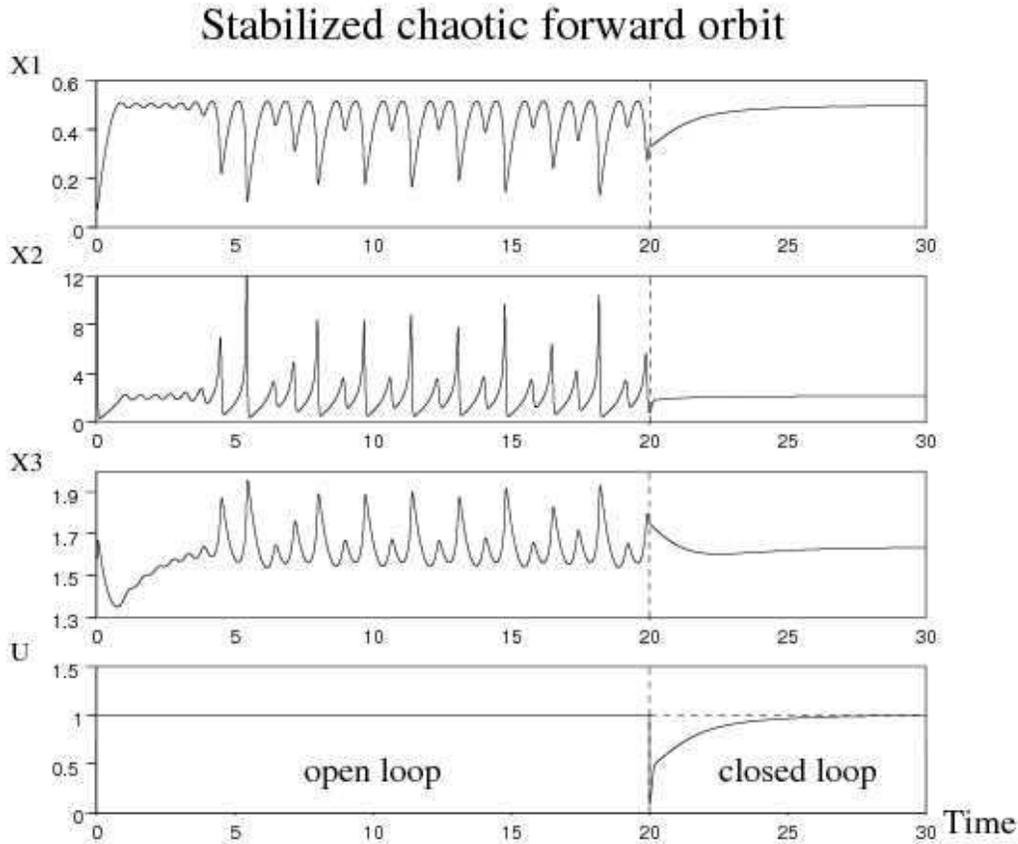}\end{center}
\caption{Time varying representation of a forward trajectory of the chaotic system $(\mathcal{S}_3)$ with $u=1$ before $t=20$ and under  control law (\ref{con_cha_otic}) from $t=20$ to the end, with $\gamma\approx 1.73>\beta_m(\approx 1.71)$ computed so that asymptotically $u(.)=\gamma\psi_3(x)$ equals its open loop value.\label{fig_example_4}}
\end{figure}

As illustrated on figure \ref{fig_example_4}, we get the following result from Theorem \ref{theo_main_result}:
 \begin{prop} For all $\gamma>0$ the control law 
\beq\label{con_cha_otic} u(.)=\gamma \psi_3(x)\eeq
globally stabilizes, on \emph{int}$(\mathbb{R}^n_+)$, system $(\mathcal{S}_3)$ towards $x^\star_\gamma$, the (sole) solution of $\gamma f_3(x)+c_3=0$.
\end{prop}

\section{Discussion}

In this paper we have considered a rather broad class of dynamical systems encompassing several possible dynamical behaviours. These models can especially represent biological systems exhibiting multistability, periodic solutions or even chaotic behaviour. In practice, the key point for stabilising these systems is that the adequate scalar function $\psi(x)$, whose exact analytical expression does not need to be known, is measured. In the example of the anaerobic digester that can have 3 steady states, the biogaz outflow rate $\psi(x)$, is easily on-line  measured though its exact analytical expression depending on bacterial activity is rather uncertain.

The proposed control law, although very simple, has proven its efficiency to globally asymptotically stabilize even chaotic systems which is known to be a rather tricky task when a part of the model ({\it i.e.} $\psi(x)$) is not perfectly known. Moreover the control law naturally fulfils the input non negativity constraints which is classical {\it e.g.} for biological systems.

The controller relevance was demonstrated in a real case for an anaerobic wastewater treatment plant (Mailleret et al., 2004). Despite its complexity (this ecosystem is based on more than 140 interacting bacterial species \citep{Delbes2001}) and the great uncertainties that characterize the biological models, the controller proved to efficiently stabilize the system. This approach was also applied to a photobioreactor where microalgae consuming nitrate and light were growing \citep{Mailleret2005}.

Another important issue when dealing with complex system is the difficulty to on-line measure the variable that must be regulated. The proposed control approach has thus the advantage of proposing a regulation scheme that does not require this measurement, which can be very difficult to carry out for biological systems. Moreover, and as it has been demonstrated with the real implementation, it is rather robust and stabilizes the system even if the model parameters are not perfectly known. The price to pay is however the possible lack of accuracy, since the solution of $\gamma f(x)+c=0$ can differ from the expected set point in case of parametric uncertainty. In such a case the idea consists in adding an integrator to correct the bias. This adaptive approach was proposed in (Mailleret et al., 2004) for the specific case of anaerobic wastewater treatment plant and the convergence of the controller was demonstrated. Of course it requires more measurement capability on the system, but we assume that the adaptation can be slower and thus needs a lower sampling frequency. The convergence proof of the adaptive controller in the general case 
is now the next challenge that we will tackle.

{\bf Acknowledgements:}{ The authors would like to thank V. Lemesle (COM - University of Marseille, France) and F. Grognard %\linebreak
(COMORE - INRIA,  France) for fruitful discussions on the subject.}

\bibliography{./biblio_tot}

\begin{thebibliography}{10}

\bibitem{posta2003}
L.~Benvenuti, A.~De~Santis, and L.~Farina, editors.
\newblock {\em Positive Systems: Theory and Applications}, volume 294 of {\em
  Lecture Notes in Control and Information Sciences}.
\newblock Springer, 2003.

\bibitem{De2001a}
P.~De~Leenheer and D.~Aeyels.
\newblock Stability properties of equilibria of classes of cooperative systems.
\newblock {\em IEEE Transactions on Automatic Control}, 46(12):1996--2001,
  December 2001.

\bibitem{Delbes2001}
C.~Delbes, R.~Moletta, and J.~J. Godon.
\newblock Bacterial and archaeal 16s rdna and 16s rrna dynamics during an
  acetate crisis in an anaerobic digestor ecosystem.
\newblock {\em FEMS Microbiology Ecology}, 35(1):19--26, March 2001.

\bibitem{Farina2000}
L.~Farina and S.~Rinaldi.
\newblock {\em Positive linear systems, theory and applications}.
\newblock John Wiley and Sons, 2000.

\bibitem{Frobenius1908}
G.~Frobenius.
\newblock Uber {M}atrizen aus positiven {E}lementen.
\newblock {\em Sitzungsberichte, {K}onigl. {P}reussichen {A}kad. {W}iss.},
  8:471--476, 1908.

\bibitem{Goodwin1965}
B.~C. Goodwin.
\newblock Oscillatory behavior in enzymatic control processes.
\newblock {\em Advances in Enzyme Regulation}, 3:425--438, 1965.

\bibitem{Hastings1977}
S.~Hastings, J.~Tyson, and D.~Webster.
\newblock Existence of periodic solutions for negative feedback cellular
  control systems.
\newblock {\em Journal of Differential Equations}, 25:39--64, 1977.

\bibitem{Hirsch1982}
M.~W. Hirsch.
\newblock Systems of ordinary differential equations which are competitive or
  cooperative {I}: {L}imit sets.
\newblock {\em S{IAM} {J}. {M}ath. {A}nal.}, 13:167--179, 1982.

\bibitem{Hirsch1985}
M.~W. Hirsch.
\newblock Systems of ordinary differential equations which are competitive or
  cooperative {II}: convergence almost everywhere.
\newblock {\em S{IAM} {J}. {M}ath. {A}nal.}, 16:423--439, 1985.

\bibitem{Kurzweil1956}
J.~Kurzweil.
\newblock On the inversion of {L}japunov's second theorem on stability of
  motion.
\newblock {\em American Mathematical Society Translations, Series 2},
  24:19--77, 1956.

\bibitem{Luenberger1979}
D.~G. Luenberger.
\newblock {\em {{I}ntroduction to {D}ynamic {S}ystems. {T}heory, {M}odels and
  {A}pplications}}.
\newblock John Wiley and Sons, New York, 1979.

\bibitem{Mailleret2004}
L.~Mailleret.
\newblock {\em Stabilisation globale de syst{\`e}mes dynamiques positifs mal
  connus. {A}pplications en biologie.}
\newblock PhD thesis, Universit{\'e} de Nice Sophia Antipolis, 2004.

\bibitem{Mailleret2004a}
L.~Mailleret, O.~Bernard, and J.~P. Steyer.
\newblock Nonlinear adaptive control for bioreactors with unknown kinetics.
\newblock {\em Automatica}, 40(8):1379--1385, August 2004.

\bibitem{Mailleret2005}
L.~Mailleret, J.~L. Gouz{\'e}, and O.~Bernard.
\newblock Nonlinear control for algae growth models in the chemostat.
\newblock {\em Bioprocess and {B}iosystems {E}ngineering}, 27(5):319--327,
  October 2005.

\bibitem{Metzler1945}
L.A. Metzler.
\newblock Stability of muliple markets: the {H}icks conditions.
\newblock {\em Econometrica}, 13:277--292, 1945.

\bibitem{Petrov1992}
V.~Petrov, S.K. Scott, and K.~Showalter.
\newblock Mixed-mode oscillations in chemical systems.
\newblock {\em Journal of {C}hemical {P}hysics}, 97-9:6191--6198, 1992.

\bibitem{Piccardi2002}
C.~Piccardi and S.~Rinaldi.
\newblock Remarks on excitability, stability and sign of equilibria in
  cooperative systems.
\newblock {\em Systems \& {C}ontrol {L}etters}, 46(3):153--163, July 2002.

\bibitem{Silva-Navarro1997}
G.~Silva-Navarro and J.~Alvarez-Gallegos.
\newblock Sign and stability of equilibria in quasi-monotone positive nonlinear
  systems.
\newblock {\em I{EEE} {T}ransactions on {A}utomatic {C}ontrol}, 42(3):403--407,
  March 1997.

\bibitem{Smith1986a}
H.~L. Smith.
\newblock Cooperative systems of differential equations with concave
  nonlinearities.
\newblock {\em Nonlinear analysis, {T}heory, {M}ethods and {A}pplications},
  10:1037--1052, 1986.

\bibitem{Smith1986b}
H.~L. Smith.
\newblock On the asymptotic behavior of a class of deterministic models of
  cooperating species.
\newblock {\em S{IAM} {J}ournal on {A}pplied {M}athematics}, 46:368--375, 1986.

\bibitem{Smith1995}
H.~L. Smith.
\newblock {\em Monotone dynamical systems, an introduction to the theory of
  competitive and cooperative systems}.
\newblock Mathematical Surveys and Monographs. American mathematical society,
  1995.

\bibitem{Smith1995a}
H.~L. Smith and P.~Waltman.
\newblock {\em The theory of the chemostat: dynamics of microbial competition}.
\newblock Cambridge University Press, 1995.

\bibitem{Zeidler1985}
E.~Zeidler.
\newblock {\em Nonlinear {F}unctional {A}nalysis and its {A}pplications. {I}:
  {F}ixed-{P}oint {T}heorems}.
\newblock Springer-Verlag, 1985.

\end{thebibliography}

\end{document}